\documentclass[letterpaper,reqno,11pt,oneside]{amsart} 
\pdfoutput=1

\usepackage{amsmath,amsthm,amsfonts,amssymb}
\usepackage[extension=pdf]{hyperref}
\usepackage{enumitem,comment}
\usepackage[LGR,T1]{fontenc}
\usepackage{times}
\usepackage{soul}
\usepackage{mathtools}
\usepackage{mathrsfs}
\usepackage{cases}
\usepackage{bm}
\usepackage[dvipsnames]{xcolor}
\usepackage{microtype}
\usepackage{xspace}
\usepackage[totalheight=9.6in,totalwidth=6.15in,centering]{geometry}
\usepackage{graphics,graphicx}
\usepackage{xr}
\usepackage{ifthen}
\usepackage{cases}
\usepackage{graphics}
\usepackage{wrapfig}

\usepackage[style=alphabetic,sorting=nyt,natbib=true,maxnames=99,isbn=false,doi=false,url=false,firstinits=true,hyperref=auto,arxiv=abs,backend=bibtex]{biblatex}
\AtEveryBibitem{%
\clearlist{language}}
\DeclareFieldFormat[article,inbook,incollection,inproceedings,patent,thesis,unpublished]{title}{#1\isdot}
\renewbibmacro{in:}{%
\ifentrytype{article}{}{%
 \printtext{\bibstring{in}\intitlepunct}}} 
\defbibheading{apa}[\refname]{\section*{#1}}
\DeclareFieldFormat{sentencecase}{\MakeSentenceCase{#1}} \renewbibmacro*{title}{%
\ifthenelse{\iffieldundef{title}\AND\iffieldundef{subtitle}}{}
{\ifthenelse{\ifentrytype{article}\OR\ifentrytype{inbook}%
   \OR\ifentrytype{incollection}\OR\ifentrytype{inproceedings}%
   \OR\ifentrytype{inreference}} {\printtext[title]{%
     \printfield[sentencecase]{title}%
     \setunit{\subtitlepunct}%
     \printfield[sentencecase]{subtitle}}}%
 {\printtext[title]{%
     \printfield[titlecase]{title}%
     \setunit{\subtitlepunct}%
     \printfield[titlecase]{subtitle}}}%
 \newunit}%
\printfield{titleaddon}}

\definecolor{darkblue}{rgb}{0.13,0.13,0.39}
\hypersetup{colorlinks=true,urlcolor=darkblue,citecolor=darkblue,linkcolor=darkblue,%
  pdftitle=Integrable fluctuations in the KPZ universality class,pdfauthor={D. Remenik}}

\DeclareSymbolFont{MnSyC}{U}{MnSymbolC}{m}{n}
\DeclareFontFamily{U}{MnSymbolC}{}
\DeclareFontShape{U}{MnSymbolC}{m}{n}{
    <-6>  MnSymbolC5
   <6-7>  MnSymbolC6
   <7-8>  MnSymbolC7
   <8-9>  MnSymbolC8
   <9-10> MnSymbolC9
  <10-12> MnSymbolC10
  <12->   MnSymbolC12}{}
\DeclareMathSymbol{\lefthalfcap}{\mathop}{MnSyC}{185}

\newcommand{\I}{\mathrm{i}}
\newcommand{\pp}{\mathbb{P}}

\newcommand{\ee}{\mathbb{E}}
\newcommand{\rr}{\mathbb{R}}
\newcommand{\nn}{\mathbb{N}}
\newcommand{\zz}{\mathbb{Z}}
\newcommand{\p}{\partial}
\newcommand{\uno}[1]{\mathbf{1}_{#1}}
\newcommand{\ep}{\varepsilon}

\newcommand{\qqand}{\qquad\text{and}\qquad}

\renewcommand{\d}{\mathrm{d}}
\newcommand\eqdistr{\stackrel{\uptext{\scriptsize (d)}}{=}}
\let\Re\relax
\DeclareMathOperator{\Re}{Re}
\newcommand{\uptext}[1]{\text{\upshape{#1}}}

\newcommand{\ts}{\hspace{0.1em}}
\newcommand{\tts}{\hspace{0.05em}}
\newcommand{\tsm}{\hspace{-0.1em}}
\newcommand{\ttsm}{\hspace{-0.05em}}
\makeatletter
  \let\origsubsection\subsection
  \renewcommand\subsection{\@ifstar{\starsubsection}{\nostarsubsection}}
  \newcommand\nostarsubsection[1]
  {\subsectionprelude\origsubsection{#1}\subsectionpostlude}
  \newcommand\starsubsection[1]
  {\subsectionprelude\origsubsection*{#1}\subsectionpostlude}
  \newcommand\subsectionprelude{%
    \vspace{-2pt}
  }
  \newcommand\subsectionpostlude{%
  }
\makeatother

\newcommand{\fh}{\mathfrak{h}}
\newcommand{\fg}{\mathfrak{g}}
\newcommand{\ff}{\mathfrak{f}}
\newcommand{\ft}{\mathbf{t}}

\newcommand{\fx}{\mathbf{x}}
\newcommand{\fB}{\mathbf{B}}

\newcommand{\Qt}{Q^*}
\renewcommand{\P}{\chi}
\newcommand{\bP}{\bar\chi}
\newcommand{\fK}{\mathbf{K}}
\newcommand{\fI}{\mathbf{I}}
\newcommand{\fT}{\mathbf{S}}
\newcommand\epi{\operatorname{\uptext{epi}}}
\newcommand\hypo{\operatorname{\uptext{hypo}}}
\newcommand\tr{\operatorname{\uptext{tr}}}
\newcommand\Ai{\operatorname{\uptext{Ai}}}
\newcommand\UC{\operatorname{\uptext{UC}}}

\newcommand{\aip}{\mathcal{A}}

\newcommand{\cD}{\mathcal{D}}

\numberwithin{equation}{section}

\newtheorem{thm}{Theorem}[section]

\theoremstyle{definition}

\makeatletter
	\newcommand\RedeclareMathOperator{%
	  \@ifstar{\def\rmo@s{m}\rmo@redeclare}{\def\rmo@s{o}\rmo@redeclare}%
	}
	\newcommand\rmo@redeclare[2]{%
	  \begingroup \escapechar\m@ne\xdef\@gtempa{{\string#1}}\endgroup
	  \expandafter\@ifundefined\@gtempa
	     {\@latex@error{\noexpand#1undefined}\@ehc}%
	     \relax
	  \expandafter\rmo@declmathop\rmo@s{#1}{#2}}
	\newcommand\rmo@declmathop[3]{%
	  \DeclareRobustCommand{#2}{\qopname\newmcodes@#1{#3}}%
	}
	\@onlypreamble\RedeclareMathOperator
\makeatother
\RedeclareMathOperator{\det}{\mathop{\uptext{det}}}
\RedeclareMathOperator{\ker}{\mathop{\uptext{ker}}}
\RedeclareMathOperator{\exp}{\mathop{\uptext{exp}}}
\RedeclareMathOperator{\log}{\mathop{\uptext{log}}}
\RedeclareMathOperator*{\lim}{\mathop{\uptext{lim}}}
\RedeclareMathOperator*{\sup}{\mathop{\uptext{sup}}}
\RedeclareMathOperator*{\limsup}{\mathop{\uptext{lim\hspace{1pt}sup}}}
\RedeclareMathOperator*{\liminf}{\mathop{\uptext{lim\hspace{1pt}inf}}}
\RedeclareMathOperator*{\max}{\mathop{\uptext{max}}}
\RedeclareMathOperator*{\inf}{\mathop{\uptext{inf}}}
\RedeclareMathOperator*{\min}{\mathop{\uptext{min}}}
\RedeclareMathOperator*{\cos}{\mathop{\uptext{cos}}}
\RedeclareMathOperator*{\sin}{\mathop{\uptext{sin}}}
\RedeclareMathOperator*{\arg}{\mathop{\uptext{arg}}}
\RedeclareMathOperator{\Re}{\uptext{Re}}
\RedeclareMathOperator{\Im}{\uptext{Im}}

\begin{document}

\title{Integrable fluctuations in the KPZ universality class}

\author{Daniel Remenik} \address[D.~Remenik]{
  Departamento de Ingenier\'ia Matem\'atica and Centro de Modelamiento Matem\'atico (IRL-CNRS 2807)\\
  Universidad de Chile\\
  Av. Beauchef 851, Torre Norte, Piso 5\\
  Santiago\\
  Chile} \email{dremenik@dim.uchile.cl}

\begin{abstract}
The KPZ fixed point is a scaling invariant Markov process which arises as the universal scaling limit of a broad class of models of random interface growth in one dimension, the one-dimensional KPZ universality class.
In this survey we review the construction of the KPZ fixed point and some of the history that led to it, in particular through the exact solution of the totally asymmetric simple exclusion process, a special solvable model in the class.
We also explain how the construction reveals the KPZ fixed point as a stochastic integrable system, and how from this it follows that its finite dimensional distributions satisfy a classical integrable dispersive PDE, the Kadomtsev-Petviashvili (KP) equation.
\end{abstract}

\maketitle

\section{The KPZ universality class}

The subject of this survey are the universal fluctuations of a large collection of models known as the one-dimensional \emph{Kardar-Parisi-Zhang (KPZ) universality class}.
This class includes many physical and probabilistic models of one-dimensional random growth, as well as several other models, including directed polymers in a random potential, some interacting particle systems, stochastic reaction-diffusion equations, and random stirred fluids, all of which can be represented in terms of the evolution of a one-dimensional interface.
\emph{Universality} here refers to the idea that the long time, large scale fluctuations of all the models in the class share a common description, in the form of common scaling exponents and a common scaling limit, which are independent of the microscopic description of each model.

While the belief in KPZ universality originates in statistical physics, much of the progress in its understanding, and in particular in the description of the universal KPZ scaling limits, has been achieved in the mathematical literature, through the study of some particular models which present a striking degree of exact solvability, and borrowing methods from algebraic combinatorics, representation theory, mathematical physics and integrable systems.
These \emph{integrable probabilistic systems} comprise a sprawling subject, to which we cannot do justice in this article; we refer the interested reader instead to the reviews \cite{borodinPetrov,IntProbLectures} on this topic, as well as to \cite{halpin-healyTakeuchi,quastelSpohn} for more physical perspectives.
Our main focus will be to describe part of the work in the field which in recent years has led to a very complete description of the universal scaling limit of KPZ models and its connection with integrable systems and random matrix theory.

\subsection*{Two examples} 

We begin by introducing a simple model which is not in the KPZ universality class.
Suppose that blocks of unit height fall at each site of $\zz$ at rate $1$ (i.e. at the times of a rate $1$ Poisson process).
If the tower of blocks at each site grows independently of the others, the height $h(t,x)$ at time $t\geq0$ at each site $x\in\zz$ can be described, by the classical central limit theorem (applied to the Poisson distribution), as $h(t,x)\approx t+t^{1/2}\xi$ with $\xi$ a standard normal random variable.
In other words, the height grows linearly with time, with Gaussian fluctuations of size $t^{1/2}$.
But since sites in this model, sometimes called \emph{random deposition}, are independent, $h(t,x)$ presents no interesting spatial structure.
In order to obtain a more interesting one-dimensional interface, one can add a \emph{relaxation} mechanism to the model as follows: when a block falls over site $x$, it lands on either $x$ or any of its two nearest neighbors, whichever has the lowest height (choosing, say, uniformly in case of ties).
To first order, $h(t,x)$ still grows like $t$, but fluctuations are now of order $t^{1/4}$: the relaxation mechanism has the effect of smoothing the interface, which now presents non-trivial correlations on a spatial scale of order $t^{1/2}$ (see \cite{barabasiStanley} for a discussion).
This model belongs to what physicists call the \emph{Edwards-Wilkinson universality class} \cite{EW}, which still has Gaussian fluctuations.

\begin{figure}
\centering
\includegraphics[width=0.87\textwidth]{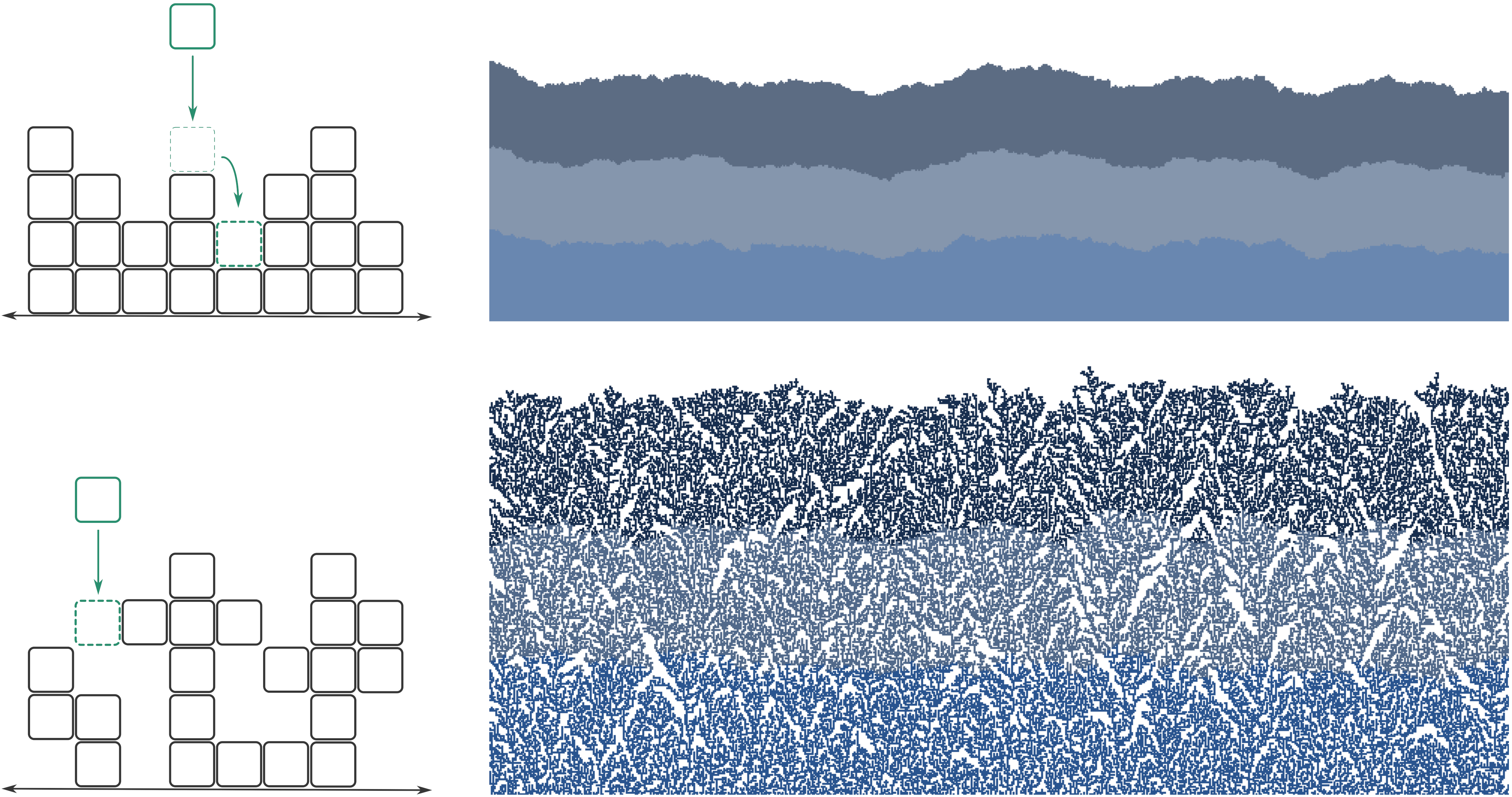}
\caption{Random deposition with relaxation (top) and ballistic deposition (bottom), simulations on the right (different shadings depict snapshots at different times).}
 \label{fig:ballistic}
\end{figure}

A very different picture arises in the \emph{ballistic deposition} model, first introduced in \cite{vold} as a model for colloidal aggregates.
In this model the falling blocks are sticky, and they attach to the side of the first neighboring block they come in contact with, see Fig. \ref{fig:ballistic}.
The interface $h(t,x)$, defined as the location of the highest block above $x$ at time $t$, now has overhangs.
\cite{seppalainenBallistic} proved that the height still grows linearly; growth has to be faster than for the previous models (as the aggregate grows, it is left with holes inside), but the exact rate remains unknown.
Fluctuations, on the other hand, are expected in this case to be of order $t^{1/3}$: the ballistic mechanism produces a rougher interface than random deposition with relaxation.
In the same way, the lateral growth of the interface makes for a longer range of spatial correlations, expected in this case to be of order $t^{2/3}$.
These two scaling exponents are one of the hallmarks of the KPZ universality class, but for this model they remain out of reach of rigorous analysis.

Ballistic deposition is representative both of the main features of models in the class and of the difficulty in analyzing them.
This is why progress in the field has had to take place mostly through the analysis of some specific models with a very special structure.
In the next section we will describe one of the main examples among these special models, TASEP.
Now we introduce the model that gives the KPZ universality class its name: the (one-dimensional) \emph{Kardar-Parisi-Zhang equation}, which is the non-linear stochastic PDE 
\begin{equation}\label{eq:KPZ}
\partial_t h = \lambda(\partial_xh)^2  + \nu \partial_x^2h + \sigma \xi,
\end{equation}
where $\xi$ is space-time white noise while $\lambda$, $\nu$ and $\sigma$ are physical parameters.
This equation was introduced in 1986 \cite{kpz} by the physicists Kardar, Parisi and Zhang, and was conceived as the simplest (and has become the canonical) continuum model for random interface growth which incorporates the physical features of models such as ballistic deposition.
Using physical arguments (based on non-rigorous dynamic renormalization group methods), \cite{forsterNelsonStephen} had predicted that the closely related stochastic Burgers equation (essentially the equation satisfied by $\p_xh(t,x)$, which can be thought of as a much simplified model of a random stirred fluid) had fluctuations of order $t^{1/3}$ with non-trivial correlations on a scale of order $t^{2/3}$. 
Using this method, \cite{kpz} then predicted that the same should hold for the KPZ equation and for a large class of models, a fortiori identified as the KPZ universality class.

The right hand side of \eqref{eq:KPZ} identifies the main elements which loosely characterize a model in the class: local dynamics, short range randomness, a smoothing mechanism (here $\p_x^2h$), and a lateral, slope-dependent growth component (which can naturally be modeled as $F(\p_xh)$ for some $F$; the $(\p_xh)^2$ term in the equation comes from keeping only the second order term in the expansion $F(u)=F(0)+F'(0)u+\frac12F''(0)u^2$, noting that the first two terms can be removed by a change of variables in the equation).
The crucial feature here is the non-linear lateral growth term (note how it becomes macroscopically apparent in the ballistic deposition simulation in Fig. \ref{fig:ballistic}).
In fact, setting $\lambda=0$ yields the (simpler, linear) \emph{additive stochastic heat equation}, which identifies the Edwards-Wilkinson class mentioned in the previous example.

In terms of solvability, the KPZ equation lies halfway between currently intractable models such as ballistic deposition and integrable models such as TASEP (this is separate from the delicate issue of well-posedness of \eqref{eq:KPZ}, which we will not discuss, see \cite{hairer}). 

\subsection*{The KPZ universality conjecture}

The Kardar-Parisi-Zhang paper marked the beginning of a long period of intense research interest, and has been one of the main drivers for advances the field, both in the physics and in the mathematics literatures.
Numerical simulations and experiments confirmed the KPZ scaling prediction for many different systems and later on, as results on the distribution of the fluctuations of special KPZ models were first obtained, the following picture began to emerge: if $h(t,x)$ is the height function describing the evolution of the interface associated to a model in the KPZ class, then (here $\eqdistr$ denotes equality in distribution)
\begin{equation}\label{eq:spatial}
\lim_{t\to\infty}t^{-1/3}\big(h(c_1t,c_2t^{2/3}x)-c_3t)\eqdistr\aip(x)
\end{equation}
for a universal limiting process $(\aip(x))_{x\in\rr}$ which depends only on the initial data of the model (more precisely, on the limit $\lim_{\ep\to0}\ep^{1/2}h(0,c_2\ep^{-1}x)$).
The scaling on the left hand side reflects the KPZ prediction: after subtracting a linear term ($c_3t$) which represents the first-order (deterministic) linear growth of a typical KPZ interface, we obtain a random variable which fluctuates at the order of $t^{1/3}$, so we need to multiply by $t^{-1/3}$ to obtain a meaningful limit, while non-trivial correlations for two spatial points are observed when they are at distance order $t^{2/3}$ (i.e. at shorter scales the height function at the two points looks the same as $t\to\infty$, while at longer scales they become independent), so the spatial variable $x$ has to be observed at that scale to see a non-trivial spatial process.
The constants $c_1,c_2,c_3$ are model dependent; they are used to provide a common normalization.

As we will see in Sec. \ref{sec:TASEPspecial}, the description \eqref{eq:spatial} emerged at first only partially: it was initially restricted only to one-point distributions (i.e., fixed $x$ instead of the whole spatial process $\aip(x)$) and, crucially, only to some very special choices of initial data.
Remarkably, it was realized that in those special cases the fluctuations arising in KPZ models were connected to random matrix theory, although to some extent the connection remained mysterious.

\subsection*{1:2:3 scaling and the KPZ fixed point}

On the other hand, \eqref{eq:spatial} does not provide a full description, since it loses all information about the temporal evolution of the interface.
To recover it, it is convenient to introduce a parameter $\ep>0$, rescale the variables $(t,x)$ as $(\ep^{-3/2}t,\ep^{-1}x)$ as well as the height $h$ (after subtracting the first order linear growth term) as $\ep^{1/2}h$, and take $\ep\to0$ (instead of $t\to\infty$).
This is usually referred to as the \emph{1:2:3 KPZ scaling} (reflecting the ratios of the exponents associated to the size of fluctuations, space and time).
The KPZ universality conjecture, first expressed in this form in \cite{cqrFixedPt}, then asserts that for any model in the class,
\begin{equation}\label{eq:fpgeneral}
\lim_{\ep\to0}\ep^{1/2}\big(h(c_1\ep^{-3/2}t,c_2\ep^{-1}x)-c_3\ep^{-3/2}t)\eqdistr\fh(t,x)
\end{equation}
for a universal process $(\fh(t,x))_{t\geq0,x\in\rr}$ which, again, should only depend on the initial data $\fh_0(x)\coloneqq\lim_{\ep\to0}\ep^{1/2}h(0,c_2\ep^{-1}x)$ prescribed for the model.
Taking $t=1$ in this limit recovers the spatial processes prescribed in \eqref{eq:spatial}.

The limiting process $\fh(t,x)$ appearing on the right hand side of \eqref{eq:fpgeneral}, is known as the \emph{KPZ fixed point}.
Since many of the models which the process should arise as a limit of are Markovian, one expects it to be a Markov process (taking values in a suitable space of real valued curves).
The name of the process comes from the fact that, by its definition as a  limit of 1:2:3 rescaled models, it should be invariant under such rescaling: if $\fh(t,x;\fh_0)$ denotes the KPZ fixed point with initial data $\fh(0,x)=\fh_0$, then for $\alpha>0$ one expects, writing $\fh_0^{(\alpha)}(x)=\alpha^{-1}\fh_0(\alpha^2x)$, that
\begin{equation}\label{eq:123inv}
\alpha\fh(\alpha^{-3}t,\alpha^{-2}x;\fh_0^{(\alpha)})\eqdistr\fh(t,x;\fh_0).
\end{equation}
The rough picture one should have in mind is of $\fh(t,x)$ as an attracting fixed point, under the renormalization map defined by the left hand side of \eqref{eq:123inv} with $\alpha\to0$, in some (loosely defined) space of models.
As such, one can think of this fixed point alternatively as \emph{defining} the KPZ universality class (as the family of models which lie in its domain of attraction).
In other words, if the KPZ fixed point can be constructed explicitly, then \eqref{eq:fpgeneral} can be used to turn the vague characterization of membership in the KPZ universality class described above into a concrete definition.

It is worth stressing that the KPZ fixed point should not be confused with the KPZ equation, which is just one (albeit very special) member of the class; in fact, the KPZ equation is not invariant under the KPZ 1:2:3 scaling \eqref{eq:123inv} (which sends the parameters $(\lambda,\nu,\sigma)$ to $(\lambda,\alpha\nu,\alpha^{1/2}\sigma)$).

Much (though certainly not all) of the progress in the field during the last 20 years can be understood as an effort to describe the KPZ fixed point, understand its properties, and explore its connections with objects coming from random matrix theory and integrable systems.
The purpose of this review is to describe one part of this story, which in particular lead to the construction of the KPZ fixed point and its description as a stochastic integrable system.
There is a priori no reason to believe that any of this should be possible; as we will see, what comes to our rescue is the remarkable exact solvability of some special discrete models in the class, which can be used to access the KPZ fixed point in the limit.

One of the main players in this story is the totally asymmetric simple exclusion process.
We turn to it in the next section.

\section{TASEP with special initial data}\label{sec:TASEPspecial}

The (one-dimensional, continuous time) \emph{totally asymmetric simple exclusion process (TASEP)} is an interacting particle system made out of particles at positions $\cdots <X_t(2)<X_t(1)< X_t(0)< X_t(-1)<X_t(-2)<\cdots$ on $\zz\cup\{-\infty,\infty\}$ performing totally asymmetric nearest neighbour random walks with exclusion: each particle independently attempts jumps to the neighbouring site to the right at rate $1$, the jump being allowed only if that site is unoccupied.
Placing particles at $\pm\infty$ allows for systems with a rightmost and/or leftmost particle with no change of notation (such particles play no role in the dynamics).
Since its introduction (in a more general form) in 1970 by F. Spitzer \cite{spitzer}, TASEP has become one of the basic and most heavily studied out-of-equilibrium models in probability and statistical physics.

In order to associate an interface to the TASEP particle system, we let $X^{-1}_t(u) = \min \{k \in \zz : X_t(k) \leq u\}$ and define the \emph{TASEP height function} as
\begin{equation}\label{defofh}
 h(t,x) = -2\tsm\left(X_t^{-1}(x-1) - X_0^{-1}(-1) \right) - x, \quad t\geq0,~x\in\zz.
\end{equation}
In words, this fixes $h(0,0)=0$ and constructs the height function by moving up from $x$ to $x+1$ whenever there is a particle at $x$ and down from $x$ to $x-1$ if the site is empty.
By interpolating piecewise linearly (and shifting $x$ by $1/2$, which makes no difference), we can picture $h(t,x)$ as a continuous function made out of line segments of slope $+1$ above every particle and $-1$ above every hole.
The dynamics of this height function is that every local maximum \raisebox{1pt}{\rotatebox{-45}{$\lefthalfcap$}} becomes a local minimum \raisebox{3pt}{\rotatebox{135}{$\lefthalfcap$}} at rate $1$, as in the figure; in this guise the model is sometimes known as the \emph{corner growth model} (for special initial data) or \emph{restricted solid-on-solid model}.
\begin{center}
\vskip3pt
\scalebox{0.8}{
\begingroup%
  \makeatletter%
  \providecommand\color[2][]{%
    \errmessage{(Inkscape) Color is used for the text in Inkscape, but the package 'color.sty' is not loaded}%
    \renewcommand\color[2][]{}%
  }%
  \providecommand\transparent[1]{%
    \errmessage{(Inkscape) Transparency is used (non-zero) for the text in Inkscape, but the package 'transparent.sty' is not loaded}%
    \renewcommand\transparent[1]{}%
  }%
  \providecommand\rotatebox[2]{#2}%
  \newcommand*\fsize{\dimexpr\f@size pt\relax}%
  \newcommand*\lineheight[1]{\fontsize{\fsize}{#1\fsize}\selectfont}%
  \ifx\svgwidth\undefined%
    \setlength{\unitlength}{228.84321438bp}%
    \ifx\svgscale\undefined%
      \relax%
    \else%
      \setlength{\unitlength}{\unitlength * \real{\svgscale}}%
    \fi%
  \else%
    \setlength{\unitlength}{\svgwidth}%
  \fi%
  \global\let\svgwidth\undefined%
  \global\let\svgscale\undefined%
  \makeatother%
  \begin{picture}(1,0.27003697)%
    \lineheight{1}%
    \setlength\tabcolsep{0pt}%
    \put(0,0){\includegraphics[width=\unitlength,page=1]{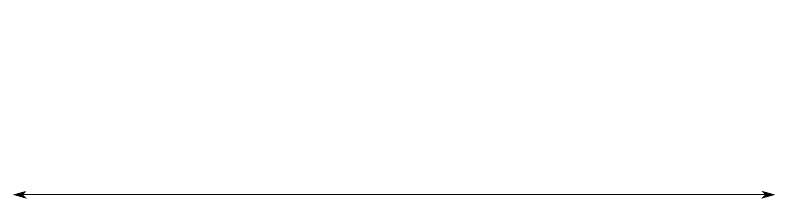}}%
    \put(0.45249215,0.22734037){\color[rgb]{0,0,0}\makebox(0,0)[lt]{\lineheight{0}\smash{\begin{tabular}[t]{l}$h(t,x)$\end{tabular}}}}%
    \put(0.36,-0.04018571){\color[rgb]{0,0,0}\makebox(0,0)[lt]{\lineheight{0}\smash{\begin{tabular}[t]{l}$X_t(N)$\end{tabular}}}}%
    \put(0.575,-0.04018571){\color[rgb]{0,0,0}\makebox(0,0)[lt]{\lineheight{0}\smash{\begin{tabular}[t]{l}$\dotsm$\end{tabular}}}}%
    \put(0.66,-0.04018571){\color[rgb]{0,0,0}\makebox(0,0)[lt]{\lineheight{0}\smash{\begin{tabular}[t]{l}$X_t(2)$\end{tabular}}}}%
    \put(0.775,-0.04018571){\color[rgb]{0,0,0}\makebox(0,0)[lt]{\lineheight{0}\smash{\begin{tabular}[t]{l}$X_t(1)$\end{tabular}}}}%
    \put(0,0){\includegraphics[width=\unitlength,page=2]{tasep.pdf}}%
  \end{picture}%
\endgroup%
}
\vskip12pt
\end{center}
To see how the features of a KPZ model described in the introduction arise in TASEP, think of writing the evolution of the height function as a stochastic equation involving a family of independent Poisson processes at each site.
Such an equation can be rewritten roughly as $\d h(t,x)=-2\uno{\text{\raisebox{1pt}{\rotatebox{-45}{$\lefthalfcap$}}}}\d t+\d M_t(x)$ with $M_t$ a martingale which provides the random forcing, and where the drift term (which says that the height function goes down by two at rate $1$ at sites where we see a local maximum) contains the smoothing and lateral growth mechanisms, as can be seen by rewriting it as
\[-2\uno{\text{\raisebox{1pt}{\rotatebox{-45}{$\lefthalfcap$}}}}=\tfrac12\left[( \nabla^-h)( \nabla^+h) ~-~ 1 ~+~ \tfrac12\nabla^+\nabla^-h\right]\]
with $\nabla^\pm$ the forward/backward discrete difference operators
\begin{equation}\label{eq:discr-diff}
\nabla^+f(x)=f(x+1)-f(x),\qquad\nabla^-f(x)=f(x)-f(x-1)
\end{equation}
(in fact, using this decomposition it can be shown that if particles now jump to the right at rate $p\in[0,1]$ and to the left at rate $q=1-p$, then the associated height function converges, in the \emph{weakly asymmetric limit} corresponding to $p-q=\ep^{1/2}$ with $\ep\to0$ under diffusive scaling, to a solution of the KPZ equation \eqref{eq:KPZ}, see \cite{berGiaco}). 

\subsection*{Wedge initial data}

The simplest possible choice of (infinite) TASEP initial condition is the \emph{packed}, or \emph{step}, initial data where particles are initially placed at every negative integer site (i.e. $X_0(i)=-i$, $i\geq1$).
For the TASEP height function it translates into the \emph{wedge} initial condition $h(0,x)=-|x|$.
The (essentially) first result about the limiting fluctuations for a KPZ model was proved by Johansson in 1999 (to be more precise, a version of this result was proved a couple of months earlier by Baik, Deift and Johansson in their seminal paper \cite{baikDeiftJohansson} for Poissonian last passage percolation):

\begin{thm}[\cite{johanssonShape}]\label{thm:TASEP-GUE}
For the TASEP height function with wedge initial data $h(0,x)=-|x|$, one has
\begin{equation}\label{eq:toGUE}
\lim_{t\to\infty}\pp\tts\bigg(\frac{h(2t,2t^{2/3}x)+t}{t^{1/3}}\leq r\bigg)=F_{\uptext{GUE}}(r-x^2),
\end{equation}
where $F_{\uptext{GUE}}$ is the Tracy-Widom GUE distribution \cite{tracyWidom}.
\end{thm}

The parabolic shift appearing on the right hand side of \eqref{eq:toGUE} reflects the curvature of the (deterministic, first order) hydrodynamic limit for the model in this case, which states \cite{rost} that $\lim_{\kappa\to\infty}\kappa^{-1}h(\kappa t,\kappa x)=-t+x^2/2t$ for $|x|\leq t$.
What is remarkable in this result, and was in fact very surprising, is the nature of the distribution of the limiting fluctuations: they coincide with the asymptotic fluctuations of the largest eigenvalue of a matrix from the \emph{Gaussian Unitary Ensemble (GUE)}, i.e. an Hermitian random matrix with properly scaled (complex) Gaussian entries.

In terms of the KPZ universality conjecture \eqref{eq:fpgeneral}, this result can be reinterpreted as the first computation of a marginal of the KPZ fixed point $\fh(t,x)$, for $t=1$, fixed $x\in\rr$ and initial data $\fh(0,x)=0$ for $x=0$ and $-\infty$ everywhere else.
We will denote this choice of initial data (which may look singular, but is natural and as we will see is in fact the simplest possible initial condition for the KPZ fixed point) by $\mathfrak{d}_0$; it is known as a \emph{narrow wedge}, as it arises from wedge initial data becoming increasingly narrower in the scaling limit.

The proof of Thm. \ref{thm:TASEP-GUE} in \cite{johanssonShape} is based on the analysis of TASEP as a determinantal point process, using as a basic tool the Robinson-Schensted-Knuth (RSK) correspondence from algebraic combinatorics.
The eigenvalues of an $N\times N$ GUE random matrix are also determinantal, and in fact it was later understood \cite{Anisotropic,warren} (see also \cite{oconnelYorNonColliding}) that using these and related tools, these eigenvalues and the location of the first $N$ TASEP particles (with step initial data) can be realized as projections of a larger process.
We will describe shortly a different proof.

The next step was to extend Thm. \ref{thm:TASEP-GUE} to the full (fixed time) spatial process:

\begin{thm}\label{thm:Airy2}
The rescaled TASEP height function $t^{-1/3}(h(2t,2t^{2/3}x)+t)$ with wedge initial data converges in distribution as $t\to\infty$, uniformly in $x$ on compact sets, to $\aip_2(x)-x^2$, where $\aip_2$ is the \emph{Airy$_2$ process}.
\end{thm}

The Airy$_2$ process was introduced by Pr\"ahofer and Spohn in 2001 \cite{prahoferSpohn} as the limit (at the level of finite dimensional distributions) of the closely related polynuclear growth (PNG) model;
a version of the result quoted here was proved by Johansson \cite{johansson} in 2003 (for a related discrete time model, see the coming discussion about LPP).
Comparing with Thm. \ref{thm:TASEP-GUE} we see that $\aip_2(x)$ has to be stationary, with Tracy-Widom GUE marginals at each $x$.
The process itself
\begin{figure}
\centering
\hskip4pt\includegraphics[width=0.5\textwidth]{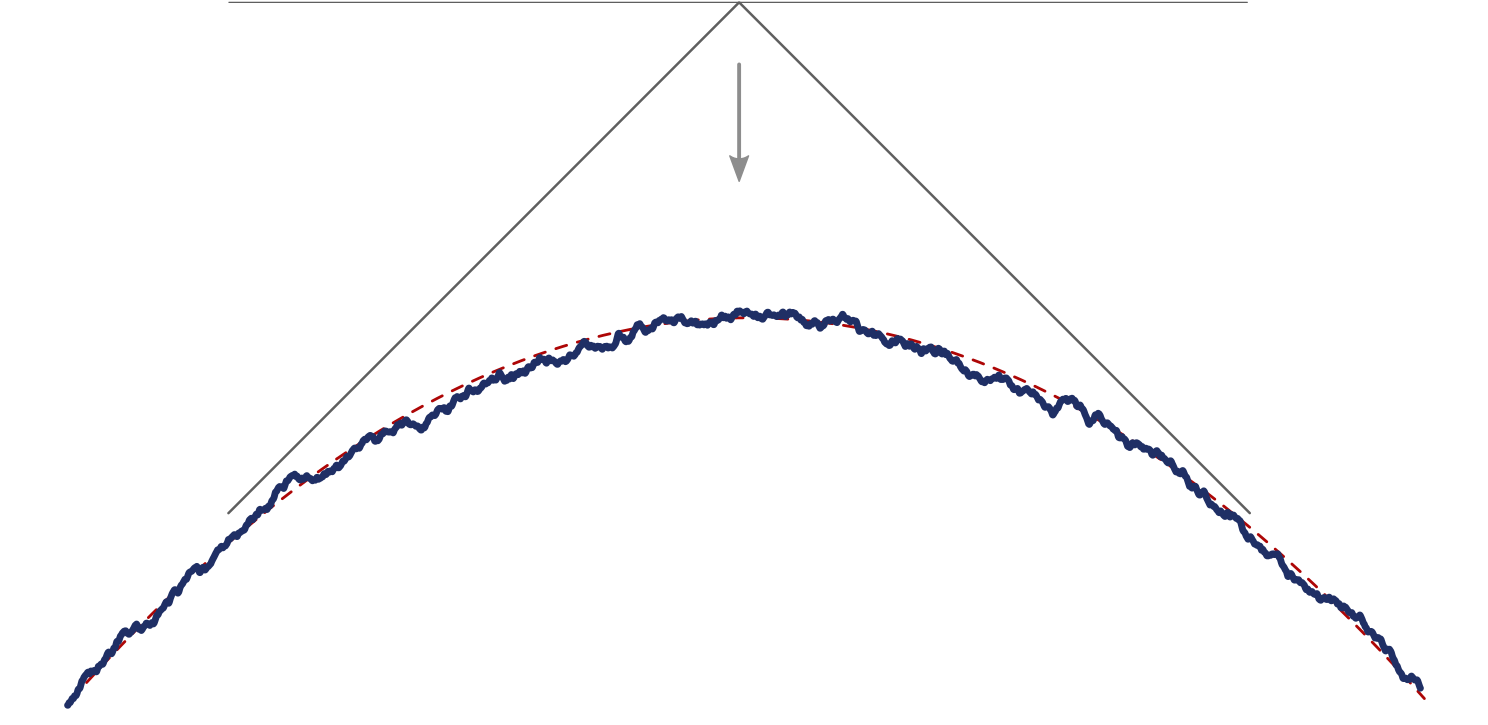}
\vskip-8pt
\caption{A simulation of the KPZ fixed point with narrow wedge initial data $\fh(1,x;\mathfrak{d}_0)\eqdistr\aip_2(x)-x^2$ as the limit of the TASEP height function}\label{fig:TASEP}
\end{figure}
is in fact closely related to random matrices.
In particular, it arises as the scaling limit of the top path of GUE Dyson Brownian motion, the eigenvalue process associated to a GUE matrix whose entries evolve as independent (complex) Brownian motions.
The process is defined through its finite dimensional distributions, which are given by a Fredholm determinant formula, see \eqref{eq:airy2} below.

In terms of the KPZ fixed point, Thm. \ref{thm:Airy2} now tells us that, as a process in $x$,
\begin{equation}
\fh(1,x;\mathfrak{d}_0)\eqdistr\aip_2(x)-x^2\label{eq:fp-a2}
\end{equation}
(where $\mathfrak{d}_0$ is the narrow wedge initial data introduced after Thm. \ref{thm:TASEP-GUE}).
See Fig. \ref{fig:TASEP}.

\subsection*{Last passage percolation}

We make a brief detour now to introduce another model in the KPZ universality class, \emph{last passage percolation (LPP)}.
We focus on the discrete case; there are similar models in other settings (e.g. Poisson LPP in the continuous case, Brownian LPP in the semi-discrete case).
Consider a family $\{w_{i,j}\}_{i,j\in\zz}$ of i.i.d. random variables and define the \emph{point-to-point last passage time}
\[L[(m_1,n_1)\to(m_2,n_2)]=\max_{\pi\in\Pi:(m_1,n_1)\to(m_2,n_2)}
\sum_i w_{\pi_i},\]
where the max is taken over the set of all paths connecting $(m_1,n_1)$ to $(m_2,n_2)$ which take unit steps up or right; if this set is empty we take the max to be $-\infty$.
We can associate a growing cluster to this model by considering the set of points with passage times less than $t$, i.e. 
\[\mathcal{C}_0(t)=\{(m,n)\in\zz^2\!:L[(0,0)\to(m,n)]\leq t\}.\]
If the $w_{i,j}$'s are exponentials with parameter $1$ then the model can be mapped to TASEP with step initial condition by intepreting each $w_{i,j}$ as the waiting time that it takes particle $j$ to jump from site $i-j$ to site $i-j+1$ (counted from the instant when that jump first becomes possible); in fact, it is not too hard to check that, after a rotation by $-3\pi/4$ (and a slight shift), the boundary of $\mathcal{C}_0(t)$ encodes the TASEP height function $h(t,\cdot)$.
In view of Thm. \ref{thm:Airy2}, the LPP fluctuations are thus governed by the Airy$_2$ process.
Johansson's result in \cite{johansson} was for the case when the $w_{i,j}$'s have a geometric distribution, which maps to a discrete time version of TASEP; it gives
\[c_1N^{-1/3}(L[(0,0)\to(N+c_2N^{2/3}x,N-c_2N^{2/3}x)]-c_3N)\xrightarrow[N\to\infty]{}\aip_2(x)-x^2\]
(note that, as stated, this LPP result is not quite equivalent to Thm. \ref{thm:Airy2} even if one lets the weights be exponential).
By universality one expects the same to hold for general choices of weights, but the problem is completely open.

To map LPP to TASEP with general initial data one can consider paths which start from any point in a given curve instead of just at the origin.
Another, perhaps more natural way of changing the LPP initial data, is to let paths start at any point in the anti-diagonal line $\{(
\ell,-\ell)\}_{\ell\in\zz}$ and add an extra reward $g(\ell)$ (the \emph{boundary condition}) there, i.e. to set
\begin{equation}
L_g[(m,n)]=\sup_{\ell\in\zz}\big(L[(\ell,-\ell)\to(m,n)]+g(\ell)).\label{eq:LPPbdry}
\end{equation}
In the scaling limit, $g$ will now become the initial data for the KPZ fixed point.

\subsection*{Periodic initial data} 

Coming back to TASEP, the next case which could be solved corresponds to \emph{periodic} initial data  $X_0(i)=-2i$, $i\in\zz$, which at the level of the TASEP height function translates into the (asymptotically) \emph{flat} initial condition of the form \raisebox{2pt}{\rotatebox{-45}{$\lefthalfcap$}}\tsm\ttsm\raisebox{2pt}{\rotatebox{-45}{$\lefthalfcap$}}\tsm\ttsm\raisebox{2pt}{\rotatebox{-45}{$\lefthalfcap$}}\tsm\ttsm\raisebox{2pt}{\rotatebox{-45}{$\lefthalfcap$}}.
This case was first solved for Poissonian LPP with zero boundary condition \cite{baikRains-symm}, which translated into the context of TASEP suggested that for periodic initial data,
\begin{equation}\label{eq:toGOE}
\lim_{t\to\infty}\pp\tts\bigg(\frac{h(2t,2t^{2/3}x)+t}{t^{1/3}}\leq r\bigg)=F_{\uptext{GOE}}(4^{1/3}r),
\end{equation}
where $F_{\uptext{GOE}}$ is the Tracy-Widom GOE distribution \cite{tracyWidom2}, the analog of $F_{\uptext{GUE}}$ for the \emph{Gaussian Orthogonal Ensemble (GOE)}, i.e. symmetric random matrices with properly scaled (real) Gaussian entries.
This was later confirmed, and extended to the full spatial process:

\begin{thm}[\cite{sasamoto,borFerPrahSasam}]\label{thm:Airy1}
The rescaled TASEP height function $t^{-1/3}(h(2t,2t^{2/3}x)+t)$ with flat initial data converges in distribution as $t\to\infty$ to the \emph{Airy$_1$ process} $\aip_1(x)$.
\end{thm}

The Airy$_1$ process is the analog of the Airy$_2$ process for flat initial data.
It is stationary, and has Tracy-Widom GOE marginals.
In terms of the KPZ fixed point, Thm. \ref{thm:Airy1} now tells us that, as a process in $x$
\begin{equation}
\fh(1,x;0)\eqdistr\aip_1(x).\label{eq:fp-a1}
\end{equation}

\subsection*{Transition probabilities}

We turn now to a sketch of the proof of Thms. \ref{thm:Airy2} and \ref{thm:Airy1}.
It is based on Sch\"utz's 1997 solution \cite{schutz} of TASEP with $N$ particles, which shows that the transition probabilities of $(X_t(1),\dotsc,X_t(N))$ have the following determinantal form:
\begin{equation}\label{eq:Green}
\pp_{X_0}(X_t(1)= x_1,\ldots,X_t(N)=x_N)=\det(F_{j-i}(t,x_{i}-X_0(j)))_{1\leq i,j\leq N},
\end{equation}
where $X_0$ in the subscript denotes the initial data of the process and where
\begin{equation}\label{eq:Fn}
F_{n}(t,x)=\frac{1}{2\pi}\oint_{\Gamma_{0,1}} \d w\,\frac{(1-w)^{-n}}{w^{x-n+1}}e^{t(w-1)},
\end{equation}
with $\Gamma_{0,1}$ any positively oriented simple loop which includes $w=0$ and $w=1$.
The derivation uses a method known in physics as the coordinate Bethe ansatz to find a solution of Kolmogorov forward (or master) equation of the process.
The key ingredients in the derivation become apparent after rewriting the functions $F_n$, $n\in\zz$, as 
\begin{equation}\label{eq:Fn2}
F_n(t,x)=(\nabla^+)^ne^{-t\nabla^-}\!\delta_0(x)
\end{equation}
where $\nabla^\pm$ are the discrete difference operators from \eqref{eq:discr-diff} (with the inverse of $\nabla^+$ defined through $(\nabla^+)^{-1}f(x)=\sum_{y>x}f(y)$) and $\delta_i(y)=\uno{y=i}$.
The operator $e^{-t\nabla^-}$ is simply the transition semigroup of a Poisson process with jumps to the left at rate $1$ (and has a kernel acting by convolution with precisely the right hand side of \eqref{eq:Fn} with $n=0$, which is where \eqref{eq:Fn2} comes from); this factor encodes the dynamics of a free (i.e., not subject to exclusion) TASEP particle.
The factor $(\nabla^+)^n$, on the other hand, encodes the exclusion restriction: very roughly put, in a situation where one particle tries to jump on top of another one, this factor produces two identical rows in the determinant obtained by using \eqref{eq:Green} on the right hand side of the Kolmogorov equation, and hence terms corresponding to those transitions will not contribute.

In principle, \eqref{eq:Fn} contains all the information one needs in order to compute a limit like \eqref{eq:fpgeneral} for TASEP, at least for initial data which has a rightmost particle $X_t(1)$.
In fact, computing the distribution of the TASEP height function at a given (finite) set of locations is equivalent to computing, for some given indices $n_1<\dotsm<n_m$,
\begin{equation}
\pp_{X_0}(X_t(n_1)>a_1,\dotsc,X_t(n_m)>a_m)\label{eq:TASEPf}
\end{equation}
and the evolution of $(X_t(i))_{i=1,\dotsc,n_m}$ is independent of the particles to their left, so we may restrict to a system with a finite number $N$ of particles.
However, \eqref{eq:Fn} is not by itself conducive to asymptotic analysis, for which we need to sum over the positions of the other $N-m$ particles and then take $N$, which is also the dimension of the determinant, to infinity.

\subsection*{Biorthogonalization}

This difficulty was overcome in \cite{sasamoto,borFerPrahSasam}, where the authors were able to show that the right hand side of \eqref{eq:Green} can be expressed as a marginal of a (signed) determinantal point process on a larger space of Gelfand-Tsetlin patterns (i.e. triangular arrays of integers with interlaced consecutive levels).
This allowed them to use techniques from random matrix theory (more precisely, a version of the Eynard-Mehta Theorem \cite{eynardMehta}) to derive an explicit 
Fredholm determinant formula for \eqref{eq:TASEPf}.
We will not describe the derivation, and content ourselves with stating a version of their result.

To do so, we need to introduce some notation.
Fix an initial condition $(X_0(n))_{n\geq1}$ for the particle system (which is right-finite, i.e. with a rightmost particle).
Define
\begin{equation}\label{eq:Q}
Q(x,y)=2^{y-x}\uno{x>y},\qquad Q^{-1}(x,y)=2^{y-x}\nabla^+(x,y)=2^{y-x}(\uno{x=y-1}-\uno{x=y});
\end{equation}
$Q$ is invertible as an operator acting on $\ell^2(\zz)$, with inverse given by $Q^{-1}$ as defined above ($\nabla^+(x,y)$ is similarly just the kernel of $\nabla^+$).
Next, for $n\geq0$ and $k<n$ let
\begin{equation}\label{eq:psin}
\Psi^n_k(x)=2^{X_0(n-k)-x}F_{-k}(t,x-X_0(n-k))=Q^{-k}e^{-\frac12t\nabla^-}\!\delta_{X_0(n-k)}(x).
\end{equation}
The powers of $2$ which we have introduced should be thought of as a convenient normalization, the crucial point being that $Q$ is the transition matrix of a random walk with strictly negative Geom$[\frac12]$ steps.
$\Psi^n_k$ can be thought of (c.f. \eqref{eq:Fn2}) as coming from applying repeatedly the forward difference operator to the Poisson weight $w_{t/2}(x)=e^{-t/2}(t/2)^{x}/x!\uno{x\geq0}$, shifted by the initial data $X_0(n-k)$.
One checks directly then, using the classical recurrence equations satisfied by the \emph{Charlier polynomials} $C_k(x,t)$ (i.e. the family of discrete orthogonal polynomials with respect to the Poisson weight $w_t(x)$) that
\[\Psi^n_k(x)=2^{X_0(n-k)-x}f_k(x+k-X_0(n-k))\quad\text{with}\quad f_k(x)=C_k(x,t)w_{t/2}(x).\]
Note that the functions $\Psi^n_k$ only depend on $n$ through a shift by the TASEP initial data.
Next for $n\geq0$ define $\{\Phi^n_k(x)\}_{k=0,\dotsc,n-1}$ as the (unique) solution of the following \emph{biorthogonalization problem}:
\begin{quote}
\vskip6pt
 \emph{%
 Given the family of \emph{shifted Charlier functions} $\{\Psi^n_k\}_{k=0,\dotsc,n-1}$, find a family of functions $\{\Phi^n_k\}_{k=0,\dotsc,n-1}$ on $\zz$ so that:
 \vskip4pt
 \begin{enumerate}[label=\uptext{(\roman*)},itemsep=0pt,leftmargin=32pt]
  	\item The two families are biorthogonal, i.e.  
   	$\sum_{x\in\zz}\Psi^n_k(x)\Phi^n_\ell(x)=\uno{k=\ell}$.
  	\item $2^{-x}\Phi^n_k(x)$ is a polynomial of degree $k$.
 \end{enumerate}
 }
\vskip6pt
\end{quote}
Finally, for a fixed vector $a\in\rr^m$ and indices $n_1<\dotsm<n_m$ let 
\begin{equation}\label{eq:defChis}
\P_a(n_j,x)=\uno{x>a_j}\qqand\bP_a(n_j,x)=\uno{x\leq a_j},
\end{equation}
which we also regard as multiplication operators (actually, projections) acting on $\ell^2(\{n_1,\dotsc,n_m\}\times\zz)$, and later also on $L^2(\{n_1,\dotsc,n_m\}\times\rr)$. 
If $a$ is a scalar we write similarly $\chi_a(x)=1-\bP_a(x)=\uno{x>a}$.

\begin{thm}[\cite{sasamoto,borFerPrahSasam}]\label{thm:bfps}
Consider TASEP with initial data $(X_0(n))_{n\geq1}$ and let $n_1,\dotsc,n_m$ be distinct positive integers.
Then for $t>0$ we have
\begin{equation}\label{eq:extKernelProbBFPS}
  \pp\!\left(X_t(n_j)>a_j,~j=1,\dotsc,m\right)=\det\!\left(I-\bar\chi_aK_t\bar\chi_a\right)_{\ell^2(\{n_1,\dotsc,n_m\}\times\zz)},
\end{equation}
where
\begin{equation}\label{eq:Kt}
K_t(n_i,x_i;n_j,x_j)=-Q^{n_j-n_i}(x_i,x_j)\uno{n_i<n_j}+Q^{n_j-n_i}K_t^{(n_i)}(x_i,x_j)
\end{equation}
with
\begin{equation}\label{eq:Ktn}
K_t^{(n)}(x,y)=\sum_{k=1}^{n}\Psi^{n}_{n-k}(x)\Phi^{n}_{n-k}(y).
\end{equation}
\end{thm}

The determinant in \eqref{eq:extKernelProbBFPS} is the \emph{Fredholm determinant}: for an integral operator $A$ acting on $L^2(X,\mu)$ with kernel $A(x,y)$,
\[\det (I+A)=\sum_{n\geq 0}\frac{1}{n!}\int_{X^{n}}\mathrm{d}\mu(x_1)\dotsm\mathrm{d}\mu(x_n)\det[A(x_{i},x_{j})]_{i,j=1}^{n}.\]
Note that the result holds for any choice of right-finite initial data $X_0$.
The point is that, if we can solve the above biorthogonalization problem for $X_0$, then we have an explicit formula for the TASEP multipoint distributions which, at least in principle, is amenable to asymptotic analysis (in fact, the size of the determinant is now fixed, and computing the scaling limit will now involve only calculating suitable limits of the kernel $K_t$).

The challenge is then to solve the above biorthogonalization problem.
One sees immediately why the choice of step/wedge initial data is the simplest in this setting.
In fact, in this case $X_0(i)=-i$, so $\Psi^n_k(x)=2^{k-n-x}f_k(x+n)$ and hence by definition the biorthogonalization problem is solved by the Charlier polynomials themselves: $\Phi^n_k(x)=c_{k}2^{x+n-k}C_k(x+n,t)$ for a suitable normalization constant $c_k$.
This leads to a relatively simple form for $K_t$ in \eqref{eq:Kt} as the \emph{(extended) Charlier kernel} (related to what in random matrix theory would be called the \emph{Charlier ensemble}), from which the TASEP limit can be extracted essentially by classical orthogonal polynomial asymptotics, leading to a kernel in terms of Airy functions (this explains the name of the Airy$_2$ process in Thm. \ref{thm:Airy2}), see \eqref{eq:airy2}.
To prove Thm. \ref{thm:Airy1} one first considers the half-periodic initial condition $X_0(i)=-2i$, $i\geq1$, in which case the biorthogonalization was solved (in 2005 \cite{sasamoto}) essentially by linear algebra (the answer \cite{borFerPrahSasam} is $\Phi^n_k(x)=c_k'2^x\sum_{\ell=0}^{n-1}\frac{t^\ell}{(2k-\ell)(\ell-1)!}\binom{2k-\ell}{k-\ell}C_\ell(x,t)$); the full periodic case is recovered by focusing on particles far to the left and taking a suitable limit.
But for more general choices of $X_0$ the method stalled for about a decade. 
(A third choice of initial data, namely a product measure, could be analyzed \cite{ferrariSpohnStat} by other methods which use crucially its stationarity for the TASEP evolution; also mixed versions of the three initial conditions could be handled, with one choice on the positive integers and another one on the negative integers).

\section{General solution of TASEP}\label{sec:TASEPgen}

As the reader can probably guess by now, what we are aiming for is to construct the KPZ fixed point as the scaling limit \eqref{eq:fpgeneral} of the TASEP height function.
Two obstacles lie in our way: we only know how to compute the limit for two special choices of initial data, and we can only do it for fixed time $t$ (we chose $t=1$ above, but other choices of fixed $t$ follow in the same way by adjusting the scaling).
To some extent, however, the two obstacles are the same: in fact, TASEP is a Markov process and we thus expect the KPZ fixed point to also be Markovian, so in order to define its temporal evolution it should be enough to characterize its transition probabilities from an arbitrary initial condition in a suitable space.

\subsection*{Biorthogonalization solution}

The general solution of the biorthogonalization problem for TASEP appeared in \cite{fixedpt}, and leads to a representation for the kernel $K_t$ from which asymptotics can be performed naturally.
Two main ingredients were used in its derivation.
The first one is the time reversal invariance satisfied by TASEP.
\begin{center}
\vskip-2pt
\scalebox{0.7}{
\begingroup%
  \makeatletter%
  \providecommand\color[2][]{%
    \errmessage{(Inkscape) Color is used for the text in Inkscape, but the package 'color.sty' is not loaded}%
    \renewcommand\color[2][]{}%
  }%
  \providecommand\transparent[1]{%
    \errmessage{(Inkscape) Transparency is used (non-zero) for the text in Inkscape, but the package 'transparent.sty' is not loaded}%
    \renewcommand\transparent[1]{}%
  }%
  \providecommand\rotatebox[2]{#2}%
  \newcommand*\fsize{\dimexpr\f@size pt\relax}%
  \newcommand*\lineheight[1]{\fontsize{\fsize}{#1\fsize}\selectfont}%
  \ifx\svgwidth\undefined%
    \setlength{\unitlength}{364.22394039bp}%
    \ifx\svgscale\undefined%
      \relax%
    \else%
      \setlength{\unitlength}{\unitlength * \real{\svgscale}}%
    \fi%
  \else%
    \setlength{\unitlength}{\svgwidth}%
  \fi%
  \global\let\svgwidth\undefined%
  \global\let\svgscale\undefined%
  \makeatother%
  \begin{picture}(1,0.14965859)%
    \lineheight{1}%
    \setlength\tabcolsep{0pt}%
    \put(0,0){\includegraphics[width=\unitlength,page=1]{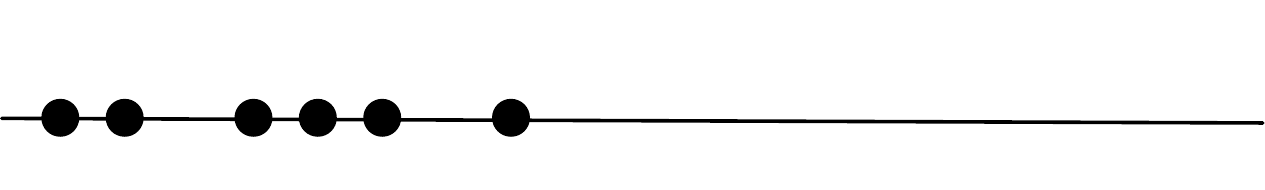}}%
    \put(0.39072497,0.01910667){\color[rgb]{0,0,0}\makebox(0,0)[lt]{\lineheight{1.25}\smash{\begin{tabular}[t]{l}$x_1$\end{tabular}}}}%
    \put(0.29144756,0.01910667){\color[rgb]{0,0,0}\makebox(0,0)[lt]{\lineheight{1.25}\smash{\begin{tabular}[t]{l}$x_2$\end{tabular}}}}%
    \put(0.0343448,0.01910667){\color[rgb]{0,0,0}\makebox(0,0)[lt]{\lineheight{1.25}\smash{\begin{tabular}[t]{l}$x_n$\end{tabular}}}}%
    \put(0,0){\includegraphics[width=\unitlength,page=2]{reversalInvarianceTASEP.pdf}}%
    \put(0.55705049,0.00748262){\color[rgb]{0,0,0}\makebox(0,0)[lt]{\lineheight{1.25}\smash{\begin{tabular}[t]{l}$a$\end{tabular}}}}%
    \put(0,0){\includegraphics[width=\unitlength,page=3]{reversalInvarianceTASEP.pdf}}%
  \end{picture}%
\endgroup%
}
\end{center}
\vskip-8pt
Suppose we start with particles at locations $x_1>\dotsm>x_n$.
By the exclusion condition, the probability that $X_t(n)>a$  is the same as the probability that $X_t(i)>a+n-i$ for each $i=1,\dotsc,n$.
But, by symmetry, this is the same as starting TASEP at $(a+1,\dotsc,a+n)$, running it backwards, and computing the probability that $X_t(i)\leq x_{n+1-i}$ for each $i$.
Using now simple reflection and shift invariance properties of the TASEP dynamics, we deduce that
\begin{equation}\label{eq:timerev}
\pp_{(x_1,\dotsc,x_n)}\big(X_t(n)>a\big)=\pp_{(-1,\dotsc,-n)}\big(X_t(i)>a-x_{n+1-i}-1,~i=1,\dotsc,n\big).
\end{equation}
We have thus turned the one-point distribution of TASEP with arbitrary initial data into the multipoint distribution of TASEP with step initial data, which as we explained in the last section can be computed explicitly.

The second ingredient is a \emph{path integral} version of the extended kernel formula \eqref{eq:extKernelProbBFPS}, which reads as follows (recall the definition of $\P_a$ for scalar $a$ after \eqref{eq:defChis}):
\begin{multline}
\pp\!\left(X_t(n_j)>a_j,~j=1,\dotsc,m\right)\\
=\det\!\big(I-K^{(n_m)}_{t}(I-Q^{n_1-n_m}\P_{a_1}Q^{n_2-n_1}\P_{a_2}\dotsm Q^{n_m-n_{m-1}}\P_{a_m})\big)_{\ell^2(\zz)},\label{eq:path-int-kernel-TASEPgem}
\end{multline}
where $K^{(n)}_t=K_t(n,\cdot;n,\cdot)$.
A formula of this type was first derived in \cite{prahoferSpohn} for the Airy$_2$ process and later extended to the Airy$_1$ process in \cite{quastelRemAiry1} and to a very wide class of processes in \cite{bcr}.
To see how this helps, observe that the factor $\P_{a_1}Q^{n_2-n_1}\P_{a_2}\dotsm Q^{n_m-n_{m-1}}\P_{a_m}(x,y)$ inside the determinant in \eqref{eq:path-int-kernel-TASEPgem} is nothing but the probability that a random walk with geometric steps goes from $x$ at time $n_1$ to $y$ at time $n_m$, staying above $a_1$ at time $n_1$, above $a_{2}$ at time $n_2$, etcetera.
On the other hand, through \eqref{eq:timerev} this formula computes the distribution of $X_t(n)$ with general initial data.
From this we see that the $\Phi^n_k$'s should be related to the probability that the geometric random walk hits the curve prescribed by the $a_i$'s, which together with the fact that the factor $K_t^{(n)}$ appearing in \eqref{eq:path-int-kernel-TASEPgem} is explicit (it is the one-point kernel for step initial data) allows one to try to guess the form of the functions.
Thm. \ref{thm:bfps} is then set up perfectly, because one can simply check (in fact, in just a few lines) that the guess gives the right answer, which is as follows: $\Phi^n_k(x)=(e^{\frac{t}{2}\nabla^-})^*h^n_k(0,x)$, where $h^n_k(\ell, x)$ is the unique solution to the initial--boundary value problem for the discrete backwards heat equation
\begin{subnumcases}{\label{bhe}}
(\Qt)^{-1}h^n_k(\ell,x)=h^n_k(\ell+1,x) &  $\ell<k,\,x \in \zz$,\label{bhe1}\\ 
h^n_k(k,x)=2^{x-X_0(n-k)} & $x \in \zz$,\label{bhe2}\\ 
h^n_k(\ell,X_0(n-\ell))= 0 & $\ell<k$.\label{bhe3}
\end{subnumcases}
For $x<X_0(n-k)$, $h^n_k(0,x)$ is simply the probability, starting from $x$, that the (reversed) random walk first goes above the curve $(X_0(n-\ell+1))_{\ell=1,\dotsc,n}$ at time $\ell=k+1$.

\subsection*{Explicit formula}

In order to obtain a usable formula for the TASEP kernel $K_t$ we need to perform the sum in \eqref{eq:Ktn} using the solution for the $\Phi^n_k$'s.
With a bit of work, this leads to a formula which is explicitly given in terms of random walk hitting times.

In order to state it we introduce a kernel $Q^n_{\epi(X_0)}$ which is defined as follows: $Q^n_{\epi(X_0)}(x,y)$ is the probability, starting at $x$, that the random walk with transition matrix $Q$ hits the strict epigraph of (i.e the region strictly above) the curve $(X_0(\ell+1))_{\ell=0,\dotsc,n-1}$ and ends at $y$ at time $n$.
One can check that, for fixed $x$, the mapping $y\longmapsto2^{-y}Q^n_{\epi(X_0)}(x,y)$ defines a polynomial for $y\leq X_0(n)$.

\begin{thm}[\cite{fixedpt}]\label{thm:rw-kernel}
The kernel $K^{(n)}_t$ in \eqref{eq:Ktn} can be written as follows:
\begin{equation}\label{eq:Ktn2}
K_t^{(n)}=e^{\frac{t}{2}\nabla^-}\!\!\!Q^{-n}\,\overline{Q}^n_{\epi(X_0)}e^{-\frac{t}{2}\nabla^-},
\end{equation}
where $\overline{Q}^n_{\epi(X_0)}(x,y)$ equals $2^y$ times the polynomial extension from $y\leq X_0(n)$ to all $y\in\zz$ of the kernel $2^{-y}Q^n_{\epi(X_0)}(x,y)$.
\end{thm}

The polynomial extension in the formula comes from using the above representation of the functions $h^n_k(0,x)$ as hitting probabilities for $x$ below the curve; since $2^{-x}\Phi ^{n}_{k}(x)$ has to be a polynomial, one can recoverit everywhere through this extension. The operators are relatively simple, however, and the polynomial extension can be computed explicitly.

\subsection*{Other processes}

The scheme which we have described works for a more general class of particle systems with determinantal transition functions of the form \eqref{eq:Green}.
It has been applied, in particular, for PushASEP \cite{nqr-kolmogorov}, one-sided reflected Brownian motions \cite{nqr-rbm}, and several discrete time versions of TASEP \cite{caterpillars}.

\section{The KPZ fixed point}

Recapitulating, what we would like to do now is to extract the limit in \eqref{eq:fpgeneral}, with $h$ the TASEP height function, using the explicit formula  supplied by Thms. \ref{thm:bfps} and \ref{thm:rw-kernel} (this involves a simple translation from the particle system to the height function).
In view of the scaling in Thms. \ref{thm:Airy2} and \ref{thm:Airy1}, we take $c_1=c_2=2$ and $c_3=-1$.
Let us briefly sketch how the limit arises.
Consider the factor $e^{\frac{t}{2}\nabla^-}\!\!\!Q^{-n}$ in \eqref{eq:Ktn2}.
Using the scaling from \eqref{eq:fpgeneral}, it becomes approximately $e^{\ep^{-3/2}t[-\nabla^-+\frac12\log(I+2\nabla^+)]}$ (ignoring lower order terms).
After suitably scaling the variables inside the kernel, the limit is computed on the scaled lattice $\ep^{1/2}\zz$, so $\nabla^\pm\sim \ep^{1/2}$ and therefore $-\nabla^- + \tfrac12\log(I+2\nabla^+) = -\nabla^-+\nabla^+-(\nabla^+)^2+\tfrac43(\nabla^+)^3 +\mathcal{O}(\ep^2)\sim\tfrac13\ep^{3/2}\partial^3$ after a simple Taylor expansion, where $\p$ is the derivative operator.
Similarly (or by the central limit theorem) we have $Q^{\ep^{-1}x}\sim e^{x\partial^2}$.
This tells us that, as $\ep\to0$, $e^{\frac{t}{2}\nabla^-}\!\!\!Q^{-n}$ becomes
\begin{equation}\label{eq:Stx}
\fT_{t,x}\coloneqq e^{\frac13t\p^3+x\p^2}.
\end{equation}
At first sight this operator may appear to be problematic, because the heat kernel $e^{x\p^2}$ is ill-defined for $x<0$, but in fact $\fT_{t,x}$ makes sense for all $t\neq0$ as an integral operator on a suitable domain with integral kernel (here $\Ai$ is the Airy function)
\[\fT_{t,x}(u,v)=t^{-1/3}e^{-\frac{2x^3}{3t^2}-(u-v)\frac{x}{t}}\tsm\Ai(t^{-1/3}(v-u)+t^{-4/3}x^2)\]
and it satisfies the group property $\fT_{t,x}\fT_{s,y}=\fT_{t+s,x+y}$ as long as $t$, $s$ and $t+s$ are all non-zero.
The convergence to $\fT_{t,x}$ which we have sketched can be proved by using a contour integral formula (similar to \eqref{eq:Fn}) for $e^{\frac{t}{2}\nabla^-}\!\!\!Q^{-n}$.
A similar argument works for the other factor, $\overline{Q}^n_{\epi(X_0)}e^{-\frac{t}{2}\nabla^-}$ (using a similar contour integral formula).
The scaling under which we are working, on the other hand, has the effect of rescaling the random walk inside this last kernel diffusively, so in the limit the random walk hitting times become Brownian hitting times.

\subsection*{Brownian scattering operator}

The above sketch explains how the ingredients which will make up the formula for the KPZ fixed point arise.
The actual proof of the limit, in a suitably strong sense and with appropriate estimates, involves some heavy asymptotic analysis.
The final result, in its most appealing form after some post-processing, involves a kernel which we introduce next.

The natural class of initial data for our (continuum) random growth models is $\UC$, the space of upper-semicontinuous functions $\fh\!:\rr\longrightarrow[-\infty,\infty)$ satisfying $\fh(x)\leq a+b|x|$ for some $a,b>0$ and $\fh\not\equiv-\infty$, which we endow with the topology of local Hausdorff convergence (with $[-\infty,\infty)$ compactified at $-\infty$).
The linear bound which we are assuming on the initial data is not quite optimal, but it ensures that $\fh(t,x)$ is defined for all $t>0$.
Note that the narrow wedge $\mathfrak{d}_0$ is in $\UC$.
Given $\fh\in\UC$ and $\ell_1<\ell_2$, let
\[\mathbf{P}_{\ell_1,\ell_2}^{\tts\uptext{No hit}\,\fh}(u_1,u_2)\tts \d u_2=\pp_{\fB(\ell_1)=u_1}\!\left(\fB(y)>\fh(y)\text{ on }[\ell_1,\ell_2],\,\fB(\ell_2)\in\d u_2\right)\]
with $\fB$ a Brownian motion with diffusion coefficient $2$, and define $\mathbf{P}_{\ell_1,\ell_2}^{\tts\uptext{Hit}\,\fh}=\fI-\mathbf{P}_{\ell_1,\ell_2}^{\tts\uptext{No hit}\,\fh}$.
The ($t$-dependent) \emph{Brownian scattering operator} associated to $\fh$ is
\begin{equation}\label{eq:scatt}
\fK^{\hypo(\fh)}_t=\lim_{\substack{\ell_1\to-\infty\\\ell_2\to\infty}}e^{-\frac13t\p^3+\ell_1\p^2}\mathbf{P}_{\ell_1,\ell_2}^{\tts\uptext{Hit}\,\fh}\ts e^{\frac13t\p^3-\ell_2\p^2}.
\end{equation}
In words, the Brownian scattering operator computes a sort of asymptotic ``transition density''
for a Brownian motion in the whole line, killed if it does not hit $\hypo(\fh)$, the hypograph of (i.e. the region below) $\fh$.
The fact that the right hand side of \eqref{eq:scatt} makes sense is far from obvious; it was first proved (for a more restricted class of $\fh$) in \cite{flat}.
A more explicit formula can be given in terms of the operators $\fT_{t,x}$ and the law of the hitting time by a Brownian motion of $\hypo(\fh)$.

The Brownian scattering operator $\fK^{\hypo(\fh)}_t$ plays a key role in the construction of the KPZ fixed point, with the function $\fh$ as the initial data for $\fh(t,x)$.
One may wonder about whether the limit on the right hand side of \eqref{eq:scatt} contains all the necessary information about $\fh$.
It does: for any $t>0$, $\fh\longmapsto\fK^{\hypo(\fh)}_t$ is invertible, and moreover continuous, as mapping from $\UC$ to a suitable space of trace class operators (the invertibility is obtained essentially directly from \eqref{eq:onesideext2} below).

\subsection*{Definition of the KPZ fixed point}

We are finally ready to define our main object of interest: the \emph{KPZ fixed point} is the (unique) Markov process taking values in $\UC$ whose transition probabilities satisfy (here $\fh_0$ in the subscript denotes the initial condition)
\begin{equation}
\pp_{\fh_0}\!\left(\fh(t,x_1)\leq r_1,\dotsc,\fh(t,x_m)\leq r_m\right)
=\det\!\left(\fI-\chi_{r}\fK^{\hypo(\fh_0)}_{t,\uptext{ext}}\chi_{r}\right)_{L^2(\{x_1,\dotsc,x_m\}\times\rr)}\label{eq:onesideext2}
\end{equation}
for any $r=(r_1,\dotsc,r_m)\in\rr^m$, with $\fK^{\hypo(\fh_0)}_{t,\uptext{ext}}$ the \emph{extended Brownian scattering operator}
\begin{equation}\label{eq:Kt-ext}
\fK^{\hypo(\fh_0)}_{t,\uptext{ext}}(x_i,u_i;x_j,u_j)=-e^{(x_j-x_i)\p^2}\!(u_i,u_j)\uno{x_i<x_j}+e^{-x_i\p^2}\fK^{\hypo(\fh_0)}_te^{x_j\p^2}\!(u_i,u_j)
\end{equation}
and where, we recall, $\P_r$ was defined in \eqref{eq:defChis}.
Two (related) statements are implicit in this definition: that the right hand side of \eqref{eq:onesideext2} defines uniquely a probability measure on $\UC$, and that it in fact defines the transition kernel of a Markov process.
The first one holds because the events in the probability on the left hand side generate the Borel $\sigma$-algebra on $\UC$.
The second one is proved based on the fact, stated next, that $\fh(t,x)$ arises as the scaling limit of TASEP, which is Markovian, plus a compactness argument which allows one to show that the property is preserved in the limit.

\begin{thm}[\cite{fixedpt}]\label{thm:TASEPtoFP}
Define
\[\fh^{\ep}(t,x)=\ep^{1/2}[h(2\ep^{-3/2}t,2\ep^{-1}x)+2\ep^{-1}x],\]
the 1:2:3 rescaled TASEP height function.
Fix $\fh_0\in\UC$ and assume that $\fh^{\ep}(0,\cdot)\longrightarrow\fh_0$ in distribution in $\UC$.
Then $\fh^\ep(t,x)\longrightarrow\fh(t,x)$ as $\ep\to0$, in distribution in $\UC$ as a process in $t,x$, where $\fh(t,x)$ is the KPZ fixed point started at $\fh_0$, i.e. the $\UC$-valued Markov process defined through \eqref{eq:onesideext2}.
\end{thm}

\subsection*{Properties of the KPZ fixed point}

The formula for the KPZ fixed point transition probabilities \eqref{eq:onesideext2} looks perhaps too complicated to be of any use, but in fact it can be used to derive many of the conjectured properties of the KPZ fixed point (as well as some suprising ones, as we will see in the next section), including the following (some of which we state vaguely, see \cite{fixedpt} for the details):
\begin{itemize}[itemsep=2pt]
\item $\fh(t,x)$ is 1:2:3 scaling invariant, i.e. it satisfies \eqref{eq:123inv}.
\item $\fh(t,x)$ is invariant under spatial shifts, reflections and affine translations.
\item Skew time reversibility: $\pp_\fg\big(\fh(\ft,\fx)\le -\ff(\fx)\big) =\pp_\ff\big(\fh(\ft,\fx)\le -\fg(\fx)\big)$ for any $\ff,\fg\in\UC$.
\item $\fh(t,x)$ is H\"older-$\frac12-$ in $x$ and H\"older-$\frac13-$ in time.
\item Brownian invariance and ergodicity: If $\fB$ is a two-sided Brownian motion with diffusion coefficient $2$ then for any $t>0$, the process $x\longmapsto\fh(t,x;\fB)-\fh(t,0;\fB)$ has the same distribution as $\mathbf{B}$.
Moreover, for any initial condition and any fixed $t>0$, the finite dimensional distributions of $\fh(t,x)-\fh(t,0)$ are locally Brownian and, under some conditions, they converge as $t\to\infty$ to those of $\mathbf{B}$.
\end{itemize}

Formula \eqref{eq:onesideext2} can also be used to recover the Airy$_1$ and Airy$_2$ processes which had already been derived for special initial data (see \eqref{eq:fp-a2} and \eqref{eq:fp-a1}), because the Brownian hitting probabilities in \eqref{eq:scatt} are explicit in those cases.
For example, for narrow wedge initial data $\mathfrak{d}_0$ the only way to hit $\hypo(\mathfrak{d}_0)$ is for the Brownian path to pass below the origin at time $0$, so one trivially has, for $\ell_1<0<\ell_2$, $\mathbf{P}_{\ell_1,\ell_2}^{\tts\uptext{Hit}\,\mathfrak{d}_0}=e^{-\ell_1\p^2}\bP_0e^{\ell_2\p^2}$ and therefore $\fK^{\hypo(\mathfrak{d}_0)}_t=e^{-\frac{t}{3}\p^3}\bP_0e^{\frac{t}{3}\p^3}$ which, using \eqref{eq:Stx} and setting $t=1$, leads directly to the known formula for the Airy$_2$ process:
\begin{equation}\label{eq:airy2}
\pp\big(\aip_2(x_1)\leq r_1,\dotsc,\aip_2(x_m)\leq r_m)=\det\!\left(\fI-\P_{r}\fK^{\Ai}_{\uptext{ext}}\ts\P_{r}\right)_{L^2(\{x_1,\dotsc,x_m\}\times\rr)}
\end{equation}
with $\fK^{\Ai}_{\uptext{ext}}$ defined by the right hand side of \eqref{eq:Kt-ext} with $\fK^{\hypo(\fh_0)}_t(u,v)$ replaced by the \emph{Airy kernel} $\int_0^\infty\d\lambda\ts\Ai(x+\lambda)\Ai(y+\lambda)$.
For flat initial data the calculation involves the hitting probability by a Brownian motion of a straight line, which can be computed using the reflection principle.

\subsection*{Variational formula} 

An alternative description of the KPZ fixed point is through a variational (Hopf-Lax type) formula involving a non-trivial input noise called the \emph{Airy sheet} $\aip(x,y)$ (which is natural, for instance, from the point of view of LPP with boundary conditions, see \eqref{eq:LPPbdry}): for the KPZ fixed point starting from $\fh(0,x)=\fh_0(x)$,
\begin{equation}\label{eq:var}
\fh(t,x)  \eqdistr  \sup_{y\in\rr}\big\{ t^{1/3}\aip(t^{-2/3} x,t^{-2/3} y)- \tfrac1{t}(x-y)^2 + \fh_0(y)\big\}.
\end{equation}
$\aip(x,y)$ can be thought of as $\fh(1,x)$ starting from a narrow wedge at $y$ at time $0$, and it therefore involves coupling different initial conditions.
The construction of the KPZ fixed point from TASEP described above leads to the Airy sheet through subsequential limits: TASEP can be constructed with coupled initial conditions but, as far as is known, one loses access to explicit formulas, and hence the distribution of the Airy sheet is unknown.
This led to a problem in that it was unclear that \eqref{eq:var} even involved a unique object on the right hand side. 
This problem was overcome in \cite{DOV}, who constructed directly the Airy sheet (and, more generally, its space-time version, known as the \emph{directed landscape}) in terms of an LPP problem on the Airy line ensemble and showed that it is the scaling limit of Brownian LPP, putting the variational formula \eqref{eq:var} on a solid footing (\cite{nqr-rbm} confirmed that, as expected, both constructions define the same objects).

The methods used in \cite{DOV} are very different from those presented here (they use heavily, in particular, a version of the RSK correspondence), and provide an alternative approach to the study of the KPZ fixed point.
As an example, they have been used to prove a strong version of the local Brownian property mentioned above: the KPZ fixed point $\fh(t,x)$ is absolutely continuous (in $x$) with respect to a Brownian motion on compact intervals \cite{sarkarVirag}.

\subsection*{Convergence for other models}

We have constructed the KPZ fixed point as a scaling limit of TASEP.
With the description of this universal limit at hand, an important and natural problem that follows is to show that it is too the limit of other models conjectured to be in the class.
Since the methods we described in Sec. \ref{sec:TASEPgen} are applicable to a wider class of determinantal interacting particle systems, it is natural to expect that convergence can be proved for those too.
This has been done for the model of one-sided reflected Brownian motions \cite{nqr-rbm}, and can also be done for the variants of TASEP covered in \cite{caterpillars} (although in this last case the detailed arguments for the scaling limit have not yet been worked out).
The methods from \cite{DOV}, on the other hand, have been extended in \cite{DV} to show convergence of several (exactly solvable) LPP models.

Recently in \cite{sarkarQuastel,virag} the convergence was extended to the KPZ equation, asymmetric exclusion processes and Brownian last passage percolation.
These major results required new ideas, since those processes are not determinantal; however, the methods still rely on integrability to a certain extent, and the convergence to the KPZ fixed point for more general models (e.g. ballistic deposition or LPP with general weights) remains wide open.

\section{Integrability of the KPZ fixed point transition probabilities}

The description of the KPZ fixed point given in the last section leaves open the question as to whether it satisfies some sort of stochastic equation.
The variational formula \eqref{eq:var} gives a partial answer; however, the distribution of the Airy sheet, and the functional of the Airy line ensemble which it arises from, are non-explicit. 
In this sense, \eqref{eq:var} is not satisfying as a universal scaling invariant equation.
But we can do something else.

\subsection*{Stochastic integrability}

Recall the definition \eqref{eq:scatt} of the Brownian scattering operator, the main building block in the KPZ fixed point formulas.
Notice that, at least formally, we can write $\fK^{\hypo(\fh_0)}_t=e^{-\frac{t}{3}\p^3}\fK^{\hypo(\fh_0)}_0 e^{\frac{t}{3}\p^3}$.
Remarkably, the dependence of $\fK^{\hypo(\fh_0)}_t$ on $t$ is completely decoupled from the dependence on the initial data.
Moreover, the time evolution is linear: at the level of the extended Brownian scattering operator \eqref{eq:Kt-ext}, it satisfies the Lax equation
\begin{equation}
\p_t \fK^{\hypo(\fh_0)}_{t,\uptext{ext}}= [ -\tfrac13 \p^3, \fK^{\hypo(\fh_0)}_{t,\uptext{ext}}],\label{eq:fixedpt-lax}
\end{equation}
where the brackets denote the commutator, $[A,B]=AB-BA$ (in other words, the equation reads $\p_t \fK^{\hypo(\fh_0)}_{t,\uptext{ext}}(x_i,u_i;x_j,u_j)= -\tfrac13(\p_{u_i}^3+\p_{u_j}^3)\fK^{\hypo(\fh_0)}_{t,\uptext{ext}}(x_i,u_i;x_j,u_j)$).
The dynamics is thus trivial at the level of the kernels, and the KPZ fixed point finite dimensional distributions are recovered by projecting down via the Fredholm determinant \eqref{eq:onesideext2}.
This provides a representation for the temporal evolution of the KPZ fixed point under which the flow is linearized; from this perspective, one might say that this presents the KPZ fixed point as a \emph{stochastic integrable system} (c.f. \cite{deift-KdV}).
Note that, in view of Thms. \ref{thm:bfps} and \ref{thm:rw-kernel}, TASEP is integrable in the same sense; this is separate from (though, of course, not unrelated to) other aspects of TASEP's exact solvability such as the coordinate Bethe ansatz leading to \eqref{eq:Fn}, the algebraic Bethe ansatz leading to its diagonalization \cite{prolhac-spectrum}, or its relation with the Schur process (see e.g. \cite{IntProbLectures}).

Note that an equation of the same nature can be written for the dependence of the Brownian scattering operator on the spatial variables: one has
\begin{equation}
(\p_{x_1}+\dotsm+\p_{x_m})\fK^{\hypo(\fh_0)}_{t,\uptext{ext}}(x_i,u_i;x_j,u_j)=(\p_{u_j}^2-\p_{u_j}^2)\fK^{\hypo(\fh_0)}_{t,\uptext{ext}}(x_i,u_i;x_j,u_j).\label{eq:fp-dx}
\end{equation}
We have stated this identity in terms of the differential operator $\p_{x_1}+\dotsm+\p_{x_m}$ in order to write a simple formula, but this will actually be consequential below.

\subsection*{Kadomtsev-Petviashvili equation}

For fixed $\fh_0\in\UC$ and $m\in\nn$, let
\begin{equation}
F(t,x_1,\dotsc,x_m,r_1,\dotsc,r_m)=\pp_{\fh_0}\!\left(\fh(t,x_1)\leq r_1,\dotsc,\fh(t,x_m)\leq r_m\right)
\end{equation}
denote the $m$-point distribution of the KPZ fixed point.
We will see now that the stochastic integrability of $\fK^{\hypo(\fh_0)}_t$ leads to a description of $F$ in terms of a classical dispersive PDE.
Shifting the variables inside the Fredholm determinant in \eqref{eq:onesideext2}, $F$ can be written as
\begin{equation}\label{eq:FdetK}
F(t,x_1,\dotsc,x_m,r_1,\dotsc,r_m)=\det(\fI-\fK)
\end{equation}
where the determinant is on the $m$-fold direct sum of $L^2([0,\infty))$ (which we have identified with $L^2(\{x_1,\dotsc,x_m\}\times[0,\infty))$) and 
\[\fK_{ij}(u_i,u_j)=\fK^{\hypo(\fh_0)}_{t,\uptext{ext}}(u_i+r_i,u_j+r_j).\]
Next we introduce an $m\times m$ matrix-valued function $Q$ defined in terms of $\fK$ as follows:
\begin{equation}
Q(t,x_1,\dotsc,x_m,r_1,\dotsc,r_m)=(\fI-\fK)^{-1}\fK(0,0)
\end{equation}
(see \cite{KP} for the fact that the right hand side is well-defined).
Note that each entry of $Q$ depends on $t$ and each of the $x_i$'s and $r_i$'s, but we omit this from the notation.
Finally, let
\begin{equation}\label{eq:drdx}
\cD_r= \partial_{r_1}+\cdots + \partial_{r_n}, \qquad \cD_x= \partial_{x_1}+\cdots + \partial_{x_n}.
\end{equation}

\begin{thm}[\cite{KP}]\label{thm:kp}
Fix an initial condition $\fh_0\in\UC$ for the KPZ fixed point and define $F$ and $Q$ as above.
Then
\[\cD_r\log F = \tr Q\]
while 
$Q$ and its derivative $q= \cD_rQ$ solve the matrix Kadomtsev--Petviashvili (KP) equation
\begin{equation}\label{eq:matKP}
\p_tq+\tfrac12\tts\cD_rq^2+\tfrac1{12}\cD_r^3q+\tfrac14\cD_x^2Q+\tfrac12[q,\tts\cD_xQ]=0.
\end{equation}
In particular, for the one point marginals of the KPZ fixed point $F(t,x,r)=\pp_{\fh_0}(\fh(t,x)\leq r)$, $\phi=\partial_r^2 \log F$ satisfies the scalar KP-II equation
\begin{equation}\label{eq:KP-II}
\partial_t \phi + \tfrac12\partial_r \phi^2 + \tfrac1{12}\partial_r^3 \phi + \tfrac14\partial_r^{-1} \partial_x^2 \phi = 0.
\end{equation}
\end{thm}

The KP equation \eqref{eq:KP-II} was originally derived from studies of long waves in shallow water \cite{ablowitzSegur}, and plays the role of a natural two dimensional extension of the Korteweg--de Vries (KdV) equation.
In fact, when $\phi$ is independent of $x$, as is the case for flat initial data $\fh_0\equiv0$ in our setting (though other, non-deterministic choices are possible \cite{lnr}), it reduces to KdV,
\begin{equation}\label{eq:KdV}
\partial_t \phi + \tfrac12\partial_r \phi^2+ \tfrac1{12}\partial_r^3 \phi=0.
\end{equation}
KP (as well as its matrix version) is completely integrable and plays an important role in the Sato theory as the first equation in the KP hierarchy \cite{jimboMiwaDate}.
However, none of the previous physical derivations of KP seem to be related to our problem.
The proof of Thm. \ref{thm:kp}, sketched below, is essentially by algebra, and we have no physical intuition yet as to why it is true.
(As a separate note, we mention that \cite{prolhacRiemann} derived essentially concurrently a connection with superpositions of KP solitons for KPZ fluctuations with special initial data in a finite volume setting).

In the one point case \eqref{eq:KP-II} the initial data in our setting is $\phi(0,x,r) = 0$ for $r\ge \fh_0(x)$, $\phi(0,x,r)=-\infty$ for $r<\fh_0(x)$.
The formal $-\infty$ can be replaced by a suitable decay condition as $t\searrow0$ but, in any case, this type of initial data is very far from what known well-posedness schemes for KP can handle, and uniqueness for the solutions of \eqref{eq:KP-II} arising from KPZ growth remains open.
The initial data for the matrix version \eqref{eq:matKP} is more delicate, see \cite{KP}.

What should be most striking about the statement of Thm. \ref{thm:kp} is that the finite dimensional distributions satisfy a closed equation at all.
In retrospect, one realizes that a PDE for the evolution of the one point distributions follows, in the special case of narrow wedge and flat initial data, from scaling considerations and \eqref{eq:toGUE}/\eqref{eq:toGOE} (see below).
But for general initial data $\fh_0\in\UC$, this is way outside the scope of what had been expected.
And in any case, even if one knows that the one point distributions satisfy a closed equation, in general one does not necessarily expect there to be multipoint equations (moreover, one would typically hope at best that the multipoint distributions satisfy an hierarchy of equations, linking the $m$-point distribution to the $k$-point distributions for $k<m$); an exception is the multipoint distributions of the Airy$_2$ process, for which PDEs were expected to be satisfied, some answers had been derived in \cite{adlerVanMoPDEs,tracyWidom-AiryDiffEqns}.
It is also not too clear why the multipoint equation should be written in terms of derivatives with respect to the variables $r_1+\dotsm+r_m$ and $x_1+\dotsm+x_m$ in \eqref{eq:drdx}.

Let us very briefly sketch how Thm. \ref{thm:kp} is proved.
For simplicity, we restrict to the one point distribution $F(t,x,r)$ (in which case the proof amounts essentially to a rediscovery of an argument which had been employed before in an abstract setting, see e.g. \cite{poppeIP}).
Letting $\Phi(t,x,r)=\p_r\log(F(t,x,r))$ we have, using \eqref{eq:FdetK},
\[\Phi(t,x,r)=\p_r\log(\det(\fI-\fK))=-\tr\left((\fI-\fK)^{-1}\p_r\fK\right)=(\fI-\fK)^{-1}\fK(0,0);\]
the second equality is standard, while the third one follows from a simple computation of the trace of $(\fI-\fK)^{-1}\p_r\fK$ as $\int_0^\infty\d x\ts(\fI-\fK)^{-1}\p_r\fK(x,x)$ and the crucial fact that, by definition,
\begin{equation}\label{eq:dr}
\p_r\fK(u,v)=(\p_u+\p_v)\fK(u,v).
\end{equation}
We can now compute derivatives of $\Phi$ directly in terms of the derivatives of $\fK$ with respect to each parameter, using the last identity together with the evolution equations \eqref{eq:fixedpt-lax} and \eqref{eq:fp-dx}; putting them together leads, after a couple of pages of computations, to \eqref{eq:KP-II} (the main difficulty is the non-linear term in the KP equation, but this can be handled through a suitable integration by parts formula).

\subsection*{Tracy-Widom distributions}

One of the most striking aspects of the study of the KPZ universality class is its connection with random matrix theory. The most prominent instance of this connection is provided by \eqref{eq:toGUE}/\eqref{eq:toGOE}, which state that the one point distribution of the KPZ fixed point with narrow wedge and flat initial data are distributed, respectively, as Tracy-Widom GUE and GOE random variables.
A partial explanation for this in the GUE/narrow wedge case is the fact, which we mentioned in Sec. \ref{sec:TASEPspecial}, that some special KPZ models (with very particular initial conditions, e.g. step for TASEP) can be coupled to models which are naturally of random matrix type; in the GOE/flat case such a connection appears not to have been fully uncovered.
In any case, the correct picture one should have in mind seems to be of KPZ and random matrix theory as two separate (though related) domains which intersect, most prominently via the central role played by the Tracy-Widom distributions on both sides.

This bears the natural question as to what makes the Tracy-Widom distributions so special.
One difficulty is that the Tracy-Widom distributions themselves seem to lack any meaningful invariance.
In the KPZ setting, however, Thm. \ref{thm:kp} makes such an invariance apparent.
The crucial point is that the KP equation \eqref{eq:KP-II} is invariant (as it has to be, in view of \eqref{eq:123inv}) under the 1:2:3 rescaling
\begin{equation}\label{eq:KP-II-inv}
\phi(t,x,r) \mapsto \alpha^{-2} \phi(\alpha^{-3} t,\alpha^{-2} x,\alpha^{-1} r),\qquad \fh_0(x)\mapsto \alpha^{-1} \fh_0(\alpha^2 x).
\end{equation}
As we explain next, the Tracy-Widom distributions then appear in the context of the KPZ universality class as special self-similar solutions of KP (the key being that both narrow wedge and flat initial are invariant under the rescaling in \eqref{eq:KP-II-inv}).

Consider first the narrow wedge case.
From \eqref{eq:fp-a2} and the 1:2:3 scaling invariance \eqref{eq:123inv} we have $\fh(t,x)+x^2/t\eqdistr t^{1/3}\aip_2(t^{-2/3}x)$ where $\aip_2$ is the Airy$_2$ process, which is stationary.
In view of \eqref{eq:KP-II-inv}, it is then natural to look for a self-similar solution of the form $\phi^\uptext{nw} (t,x,r)=t^{-2/3} \psi^\uptext{nw} (t^{-1/3} r + t^{-4/3}x^2)$.
This turns \eqref{eq:KP-II} into the ODE
\begin{equation}\label{eq:KP-IIA}
(\psi^\uptext{nw})'''+12\psi^\uptext{nw}(\psi^\uptext{nw})'-4r(\psi^\uptext{nw})'-2\psi^\uptext{nw}=0.
\end{equation}
The transformation $\psi^\uptext{nw}=-u^2$ takes \eqref{eq:KP-IIA} into the Painlev\'e II equation:
\begin{equation}\label{eq:p2}
u''=ru+2u^3.
\end{equation}
Known tail estimates for the Tracy--Widom GUE distribution (i.e. the one-point marginal of the  $\mathrm{Airy} _{2}$ process) imply that, as $r\to -\infty $, one has $\phi^\uptext{nw}(t,x,r)\sim -(\tfrac{r}{2t} +\tfrac{x^2}{2t^2})$.
This picks out the Hastings-McLeod solution of \eqref{eq:p2}, $u(r)\sim -\Ai(r)$ as $r\to \infty$, and thus we get
\begin{equation}
F(t,x,r)= \exp\left\{ -\int_{\hat r}^\infty\!\d s\, (s-\hat r)u^2(s)\right\}=F_\uptext{GUE}(t^{-1/3} r+t^{-4/3}x^2)
\end{equation}
with $\hat r=\frac{r}{t^{1/3}}+\frac{x^2}{t^{4/3}}$, the last equality being the Painlev\'e II representation for $F_\uptext{GUE}$ famously derived by Tracy and Widom \cite{tracyWidom}.

In the flat case, $\fh_0=0$, there is no dependence on $x$, so we need to look for self-similar solutions of KdV \eqref{eq:KdV}, of the form $\phi^\uptext{fl}(t,r) = (t/4)^{-2/3} \psi^\uptext{fl}((t/4)^{-1/3} r)$ (the factor of $1/4$ is for convenience), leading to $(\psi^\uptext{fl})''' +12(\psi^\uptext{fl})'\psi^\uptext{fl}-r(\psi^\uptext{fl})' -2\psi^\uptext{fl} =0$.
Miura's transform $\psi^\uptext{fl} = \frac12 (u' -u^2)$
brings this to Painlev\'e II  \eqref{eq:p2}, with the same asymptotic behavior as $r\to-\infty$, and we recover, writing $\hat r=4^{1/3}t^{-1/3} r$ (the second equality comes now from \cite{tracyWidom2})
\begin{equation}
F(t,x,r)= \exp\left\{ -\frac12\int_{\hat r}^\infty\!\d s\,u(s)\right\}F_\uptext{GUE}(\hat r)^{1/2}=F_\uptext{GOE}(4^{1/3}t^{-1/3} r).
\end{equation}

\section{KP in special solutions of the KPZ equation}

The proof of Thm. \ref{thm:kp} which we sketched above does not make any use of the particular form of the Brownian scattering operator other than the fact that it satisfies the differential equations \eqref{eq:fixedpt-lax}, \eqref{eq:fp-dx} and \eqref{eq:dr} (together with some technical conditions).
The method thus applies in general to Fredholm determinants of kernels which satisfy the same identities.

Surprisingly, it turns out that, in the one point case, the method is applicable to some special solutions of the KPZ equation \eqref{eq:KPZ}, where for convenience we fix the scaling $\lambda=\nu=\tfrac14$ and $\sigma=1$.
The simplest case is, again, the narrow wedge solution $h_\uptext{nw}$ of \eqref{eq:KPZ}, by which what is meant (see \cite{berGiaco}) is that $h_\uptext{nw}=\log(Z)$ with $Z$ the fundamental solution of the stochastic heat equation with multiplicative noise (i.e. $\p_t Z = \tfrac14 \p^2_xZ + \xi Z$, $Z(0,x)=\delta_0(x)$): in the early 2010's a formula was obtained \cite{acq,sasamSpohn,dotsenkoNW,calabreseLeDousallNW} for the KPZ generating function
\begin{equation}\label{eq:gumbel}
G_\uptext{nw}(t,x,r) = \ee\left[\exp\!\left\{-e^{h_\uptext{nw}(t,x)+\frac{t}{12}-r}\right\}\right]
\end{equation}
which can be rewritten as $\det(\fI-\fK)$ with 
\[\fK(u,v)=\int_{-\infty}^\infty d\lambda\,t^{-2/3}\tfrac{e^{(v-u)x/t}}{1+e^{\lambda}}\Ai(t^{-1/3}(u+r-\lambda)+t^{-4/3}x^2)\Ai(t^{-1/3}(v+r-\lambda)+t^{-4/3}x^2).\]
This kernel satisfies the necessary equations, and from this ones gets \cite{KP} that, remarkably,
\begin{equation}\label{eq:KPKPZeqn}
\phi_\uptext{nw}\coloneqq\p_r^2\log G_\uptext{nw}\uptext{\quad solves the KP equation \eqref{eq:KP-II}}.
\end{equation}
Similar derivations are available for the KPZ equation with half-Brownian/spiked, two-sided Brownian and stationary initial data \cite{KP,zhang-BR}.
What is common to these cases is that there are special solvable models which converge to the KPZ equation and for which it has been possible to derive explicit formulas under these choices of initial data.
But this remains out of reach for general initial conditions (and for multipoint distributions), and at this point it is not known whether a statement such as \eqref{eq:KPKPZeqn} holds in any greater generality.

\vspace{10pt}

\noindent{\bf Acknowledgements.} 
Much of the work reported in this note comes from joint works with Konstantin Matetski and Jeremy Quastel.
The author is grateful to Joaquín Fontbona and Jeremy Quastel for comments on a draft of this article.
This work was partially supported by Centro de Modelamiento Matemático (CMM) Basal Funds AFB170001, ACE210010 and FB210005 from ANID-Chile, and by Fondecyt Grant 1201914.

\printbibliography[heading=apa]

\end{document}